\documentclass[graybox]{svmult}

\usepackage{type1cm}
\usepackage{makeidx}
\usepackage{graphicx} 
\usepackage{multicol}   
\usepackage[bottom]{footmisc}
\usepackage{newtxtext}       
\usepackage[varvw]{newtxmath} 
\usepackage{url}

\makeindex

\begin{document}

\title*{Stability analysis and optimal control of the logistic equation}
\toctitle{Stability analysis and optimal control of the logistic equation}
\titlerunning{Stability analysis and optimal control of the logistic equation}

\author{Márcia Lemos-Silva\orcidID{0000-0001-5466-0504} 
\and Sandra Vaz\orcidID{0000-0002-1507-2272} 
\and Delfim F. M. Torres\orcidID{0000-0001-8641-2505}\thanks{Corresponding author: delfim@ua.pt}}
\tocauthor{Márcia Lemos-Silva, Sandra Vaz and Delfim F. M. Torres}
\authorrunning{M. Lemos-Silva, S. Vaz, D. F. M. Torres} 

\institute{Márcia Lemos-Silva \at
Center for Research and Development in Mathematics and Applications (CIDMA),\\ 
Department of Mathematics, University of Aveiro,
3810-193 Aveiro, Portugal\\ 
\email{marcialemos@ua.pt}
\and 
Sandra Vaz \at
Center of Mathematics and Applications (CMA-UBI),\\
Department of Mathematics, University of Beira Interior,
6201-001 Covilh\~{a}, Portugal\\ 
\email{svaz@ubi.pt}
\and
Delfim F. M. Torres \at
Center for Research and Development in Mathematics and Applications (CIDMA),\\ 
Department of Mathematics, University of Aveiro, 3810-193 Aveiro, Portugal\\
\email{delfim@ua.pt}}

\maketitle


\abstract{We consider the unexploited/exploited logistic equation 
and study the stability of equilibrium points through Lyapunov functions.
Then, we apply first and second order optimality conditions for the optimal control
of the total biomass yield. Finally, we note that the time-scale 
logistic equation present in the literature lacks biological significance and
we propose a new version of a dynamic logistic equation, valid on an arbitrary
time scale, for which any trajectory of the system, beginning with a
positive initial condition, remains nonnegative.\newline
\newline\newline
\noindent {\bf MSC Classification 2020:} 34D23, 34N05, 49K05.}


\section{Introduction}

Differential equations can be used to represent the size of a population as it varies 
over time. In particular, the logistic equation was first introduced in the 19th century 
by the Belgian mathematician Verhulst \cite{logistic:history}. This equation generalizes 
the exponential growth equation by incorporating a maximum population value. Exponential 
growth predicts that, as time progresses, the population grows without bound. This is 
obviously unrealistic in a real-world setting, since several factors limit the rate of 
growth of a population: birth rate, death rate, food supply, predators, and so on. 
Many years ago, biologists have found that in many biological systems, the population 
grows until a certain steady-state population is reached. Indeed, as Verhulst wrote in 
\emph{A note on the law of population growth} \cite{quote}:

\begin{quote}
\emph{``We know that the famous Malthus showed the principle that the human population tends 
to grow in a geometric progression so as to double after a certain period of time, for example
every twenty five years. This proposition is beyond dispute if abstraction is made of the
increasing difficulty to find food [\dots]
The virtual increase of the population is therefore limited by the size and the fertility of
the country. As a result the population gets closer and closer to a steady state.''}	
\end{quote}

With this in mind, Verhulst proposed the following logistic growth model for the population 
$x$ at time $t$: 
\begin{equation}
\label{logistic}
\dot{x} = rx\left(1 - \frac{x}{k}\right), 
\end{equation}
where $r$ denotes the growth rate and $k$ the carrying capacity, i.e., the maximum population 
that the environment can sustain indefinitely. 

\begin{figure}[ht]
\centering
\includegraphics[scale=0.7]{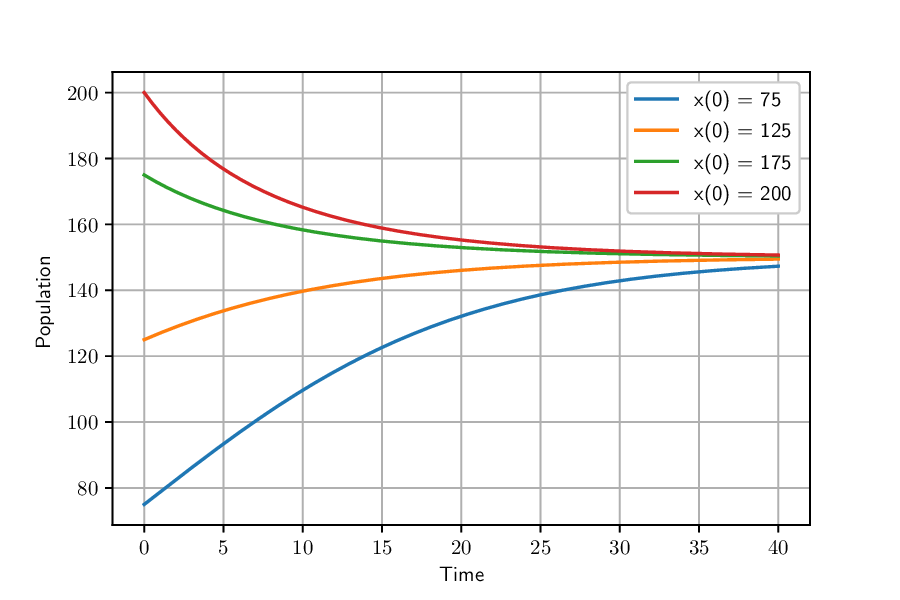}
\caption{Logistic growth \eqref{eq:log:sol}
for different values of $x(0)$ 
with $r = 0.1$ and $k = 150$.}
\label{fig1}
\end{figure}

When the population $x(t)$ is small compared to $k$, we get the approximate equation
$$\dot{x} = rx,$$
whose solution is the well-known exponential growth given by 
$$x(t) = x(0)e^{rt},$$
where $x(0)$ denotes the initial population of $x$. On the other hand, as $x(t)$ gets closer 
to $k$, $x(t)$ approaches a steady state. More than that, $\dot{x}(t)$ would actually become 
negative if $x(t)$ could exceed the value of $k$. All these behaviors can be seen in 
Figure~\ref{fig1}.

The logistic differential equation \eqref{fig1} is an autonomous differential equation, so it 
can be solved using separation of variables. This solution is well-known (see, e.g., 
\cite{logistic:history}) and is given by 
\begin{equation}
\label{eq:log:sol}
x(t) = \frac{x(0)ke^{rt}}{k + x(0)(e^{rt} - 1)},
\end{equation}
where $x(0)$ denotes the initial value for the population $x$. As one can see from this solution, 
the total population increases (or decreases) progressively from $x(0)$ to the limit $k$, which 
is reached only when $t \rightarrow +\infty$.

Despite its historical origins and apparent simplicity, the logistic equation \eqref{logistic} remains a subject of 
continued interest and it continues to be used by researchers worldwide for different applications 
(see, e.g., \cite{prey:predator,MR4708597,covid}). 

As far as we are concerned, there are few works that convey fundamental concepts on both stability
analysis and optimal control theory. Thus, the main goal of this work is to present such results in a gentle way
and apply them to the well-known logistic equation. This work is organized as follows. In Section~\ref{sec2}, 
we start by presenting some definitions on stability and Lyapunov type results to the 
analysis of global stability. We finish this section by applying these very same results to the 
standard logistic growth model and to two different exploited single-species models. In Section~\ref{sec3}, 
we state necessary conditions for a set $(x^*,u^*)$ to be optimal, known as the 
Pontryagin Maximum Principle. Moreover, we present Goh's second-order necessary conditions that are
needed when dealing with singular control problems. We end by illustrating these results on 
an optimal control problem of the logistic model. This work ends with a short section (Section~\ref{sec4}) 
explaining the relevance of this work for future work on time scales.


\section{Stability of the logistic equation}
\label{sec2}

In ecology, the word ``stability'' can be used with two different purposes: whether to describe
the lack of change in population densities or parameters, 
or to describe the persistence of a system. 
In mathematics, one can find many different concepts of stability. In this section, we will focus on 
some approaches to evaluate the local and global stability of equilibrium points.  

We begin by recalling some important definitions and results on stability for systems of ordinary 
differential equations (ODEs) of the form 
\begin{equation}
\label{def:ode}
\dot{x} = f(x),
\end{equation}
with the initial condition $x(a) = x_a$, where $a \in \mathbb{R}^+_0$, $x(t) \in \mathbb{R}^n$ 
for all $t \geq a$ and $f$ is continuously differentiable as many times as we need.

\begin{definition}[See \cite{stability}]
An equilibrium point of system \eqref{def:ode} is a vector $\tilde{x} \in \mathbb{R}^n$ that 
satisfies $f(\tilde{x}) = 0$. 
\end{definition} 

\begin{definition}[See \cite{stability}]
\label{local:def}
Let $\tilde{x}$ be an equilibrium point of system \eqref{def:ode} and $x(t)$ a 
solution of the same system satisfying the initial condition $x(a) = x_a$. Then, the equilibrium 
point $\tilde{x}$ is said to be
\begin{itemize}

\item stable if
$$
\forall \varepsilon > 0, \exists \delta > 0 :\, \rVert x_a - \tilde{x}\rVert < \delta 
\Rightarrow \lVert x(t) - \tilde{x} \rVert < \varepsilon
$$
for all $t \geq a$;

\item asymptotically stable if it is stable and  
$$\exists \delta > 0 :\, \lVert x_a - \tilde{x}\rVert < \delta 
\Rightarrow \lim_{t\rightarrow +\infty} \lVert x(t) - \tilde{x} \rVert = 0;$$
\item unstable if it is not stable. 
\end{itemize}
\end{definition}

Here, $\lVert\cdot\rVert$ denotes an arbitrary norm. 
In the present context, the Euclidean norm is frequently used.

Let us break down Definition~\ref{local:def} in simpler terms. Take $x_a$ as an initial point in a 
neighborhood of a certain equilibrium point $\tilde{x}$. We say that $\tilde{x}$ is stable if, and 
only if, the solution of system $\eqref{def:ode}$ with the initial condition $x(a) = x_a$ remains 
in a neighborhood of $\tilde{x}$. Moreover, when $\tilde{x}$ is asymptotically stable, the solution 
not only remains in a neighborhood of $\tilde{x}$ but also converges to it. On the contrary, if 
$\tilde{x}$ is unstable, then even when $x(a)$ is close to $\tilde{x}$, there is always a solution, 
say $x_u$, of system \eqref{def:ode} for which it is not possible to find a neighborhood of 
$\tilde{x}$ where $x_u$ could remain. This means that $x_u$ diverges away from $\tilde{x}$.

A stronger concept of stability is global stability. In simple terms, one has
global stability when a trajectory with an arbitrary initial point in the domain 
will remain a fixed distance from the equilibrium point;
global asymptotically stability when
any trajectory of the system 
tends to $\tilde{x}$ regardless of initial conditions 
\cite{global:stability:2013}. Thus, one can write the following definition.

\begin{definition}[See \cite{global:stability:2013}]
Let $\tilde{x}$ be an equilibrium point of system \eqref{def:ode} and $x(t)$ a 
solution of the same system satisfying the initial condition $x(a) = x_a$. 
Then, the equilibrium point $\tilde{x}$ is said to be
\begin{itemize}	
\item globally stable if for any $x(a) \in \mathbb{R}^n$							
$$
\exists \, \varepsilon > 0 \  : 
\lVert x(t) - \tilde{x} \rVert < \varepsilon
$$
for all $t \geq a$;
	
\item  globally asymptotically stable if for any 
$x(a) \in \mathbb{R}^n$							
$$
\lim_{t\rightarrow +\infty} \lVert x(t) - \tilde{x} \rVert = 0.
$$
\end{itemize}
\end{definition} 

In the world of ecological models, 
the population density of each specie must be nonnegative. This 
requires that the concept of global stability must be modified in these models. 
For example, in the asymptotic case, by definition, an ecosystem 
is globally asymptotically stable if every trajectory of the model 
that begins at a \emph{positive} state 
remains in the \emph{positive} orthant for all finite values of $t$, 
and converges to a \emph{positive}
equilibrium point as $t\rightarrow \infty$ \cite{book}.

To illustrate this, consider a model of the form 
\begin{equation}
\label{population:sys}
\dot{x}_i = x_i f_i(x_1,x_2,\dots,x_n),\quad i = 1,2,\dots,n, 
\end{equation}
where $x_i$ denotes the population density of the $i$th-specie and $f_i(x)$ are continuous 
functions in the positive orthant. Suppose \eqref{population:sys} 
has a positive equilibrium $\tilde{x}$, i.e., 
$f_i(\tilde{x}) = 0$, $i = 1,2,\dots,n$. 

For all species to persist it is necessary for model \eqref{population:sys} to have either a 
positive equilibrium or a limit cycle and all species have to be present initially. Thus, in this 
case, a natural definition for global stability of a positive equilibrium is that every solution 
of model \eqref{population:sys} that begins in the positive orthant $\mathbb{R}^n_+$ must remain in $\mathbb{R}^n_+$ and at a fixed distance from the positive equilibrium point.

The simplest way to examine stability is by evaluating the eigenvalues of the Jacobian matrix at an 
equilibrium point of the model. However, this method only determines local stability or instability. 
In contrast, real-world ecosystems experience significant perturbations and continual disturbances. 
A powerful analytical method for determining the stability of an equilibrium relative to finite 
perturbations is the direct method of Lyapunov, 
named after the mathematician 
Aleksandr Lyapunov (1820--1868) 
\cite{Liapunov:origin}.

Let $\Omega$ denote the positive region. Therefore, $\tilde{x} \in \Omega$. The following results
allow us to determine stability for a system without explicitly integrating the differential 
equation. This method is called the direct method of Lyapunov. In a physical system, a Lyapunov 
function is a generalization of the free energy of the system. If the energy is continuously 
dissipated, the state of the system will evolve and move towards the state where the energy
of the system attains local or global minimum. The corresponding state is none other than an 
equilibrium of the system.

\begin{theorem}[See \cite{book}]
\label{Lyapunov:def}
Let $V: \mathbb{R}^n \rightarrow \mathbb{R}$ be a continuously differentiable scalar function. 
If $V(x)$ satisfies the properties
\begin{itemize}
\item $V(\tilde{x}) = 0$;
\item $V(x)$ has a global minimum in $\Omega$
equal to zero at $\tilde{x}$;
\item for every $i = 1, 2, \dots, n$, $V(x) \rightarrow \infty$ as $x_i \rightarrow 0^+$ or 
$x_i \rightarrow \infty$;
\item $\dot{V}(x) = \sum_{i=1}^n \frac{\partial V}{\partial x_i} x_if_i(x) \leq 0$, 
for all $x \in \Omega$;
\end{itemize}
then the positive equilibrium $\tilde{x}$ of model \eqref{population:sys} is globally stable.
\end{theorem}

\begin{theorem}[See \cite{book}]
\label{thm:Lyapunov:stab}
Under the hypothesis of Theorem~\ref{Lyapunov:def}, the positive equilibrium 
$\tilde{x}$ of model \eqref{population:sys} 
is globally asymptotically stable in the positive orthant, if 
$$
\dot{V}(x) = \sum_{i=1}^n \frac{\partial V}{\partial x_i} x_if_i(x) < 0 \ 
\text{ for all } 
x \in \mathbb{R}^n_+ 
\text{ and } 
x \neq \tilde{x}. 
$$
\end{theorem}

\begin{definition}[See \cite{book}]
The region $\Omega$ is a finite region of stability if $\Omega$ is a finite region and all 
solutions of \eqref{population:sys} that begin in $\Omega$ remain in $\Omega$ and tend to 
$\tilde{x}$ as $t \rightarrow \infty$.  
\end{definition} 

Note that the solutions must not leave region $\Omega$ at any time $t$. Otherwise, some of 
the species may become extinct for a temporary period, which has no biological meaning.  

All results stated before are only sufficient conditions for an equilibrium to be stable. Our
inability to find a Lyapunov function $V(x)$ does not imply the instability of the equilibrium. 
For such conclusion to be made, the following theorem is needed. 

\begin{theorem}[See \cite{khalil2002nonlinear}]
\label{thm:unstable}
Let $V: \mathbb{R}^n \rightarrow \mathbb{R}$ be a continuously differentiable scalar function. 
If
\begin{itemize}
\item $V(\tilde{x}) = 0$;
\item $V(x)$ has a global minimum in $\Omega$
equal to zero at $\tilde{x}$;
\item for every $i = 1, 2, \dots, n$, $V(x) \rightarrow \infty$ as $x_i \rightarrow 0^+$ or 
$x_i \rightarrow \infty$;
\item $\dot{V}(x) = \sum_{i=1}^n \frac{\partial V}{\partial x_i} x_if_i(x) > 0$, 
for all $x \in \Omega$;
\end{itemize}
then $\tilde{x}$ is unstable.
\end{theorem}

Consider the general single-species model 
\begin{equation}
\label{single:general}
\dot{x} = xf(x).	
\end{equation}
Assume \eqref{single:general} 
has a positive equilibrium at $\tilde{x}$ and consider
\begin{equation}
\label{Lyapunov:general}
V(x) = x - \tilde{x} - \tilde{x}\ln\left(\frac{x}{\tilde{x}}\right).	
\end{equation} 
First, note that 
$$\lim_{x\rightarrow \infty} x - \tilde{x} - \tilde{x}\ln\left(\frac{x}{\tilde{x}}\right) 
= \lim_{x\rightarrow \infty} x\left(1 - \frac{\tilde{x}}{x} - 
\frac{\tilde{x}\ln(\frac{x}{\tilde{x}})}{x}\right),$$ 
and applying L'Hôpital's rule to the term $\frac{\tilde{x}\ln(\frac{x}{\tilde{x}})}{x}$, 
it is clear that
$$\lim_{x\rightarrow \infty} x\left(1 - \frac{\tilde{x}}{x} - 
\frac{\tilde{x}\ln(\frac{x}{\tilde{x}})}{x}\right) = \infty.$$
Moreover,
$$
\lim_{x\rightarrow0^+} x - \tilde{x} - \tilde{x}\ln\left(\frac{x}{\tilde{x}}\right) = 
0 - \tilde{x} -(-\infty) = \infty.
$$
It is also clear that $V(x)$ has a global minimum equal to zero at $\tilde{x}$. In fact, if we
compute the derivative of $V(x)$, we get 
$$
V'(x) = 0 \Leftrightarrow 1 - \frac{\tilde{x}}{x} = 0 \Leftrightarrow x = \tilde{x},
$$
and since $V(x) \rightarrow \infty$ as $x\rightarrow \infty$ and as $x\rightarrow 0^+$, 
it is clear that $\tilde{x}$ is a global minimizer. 
Moreover, we have that $V(\tilde{x}) = 0$. 

Furthermore, along solutions of \eqref{single:general}, we have $\dot{V}(x) = 
(x - \tilde{x})f(x)$. Thus, following 
Theorem~\ref{thm:Lyapunov:stab}, it follows that 
$\tilde{x}$ is globally asymptotically stable in the positive orthant if 
\begin{itemize}
\item for $\tilde{x} > x > 0$, $f(x) > 0$,
\end{itemize}
and
\begin{itemize}
\item for $x > \tilde{x}$, $f(x) < 0$,
\end{itemize}
for all $x \in \mathbb{R}_+$ and $x \neq \tilde{x}$.

We can state this result as a theorem. 

\begin{theorem}
\label{cond:single}
The single species model $\eqref{single:general}$ is globally asymptotically stable if
\begin{itemize}
\item it has a positive equilibrium at $\tilde{x}$;
\item for $\tilde{x} > x > 0$, $f(x) > 0$;
\item for $x > \tilde{x}$, $f(x) < 0$.   
\end{itemize}
\end{theorem}

The previous theorem allows us to state the following corollary. 

\begin{corollary}
If $\tilde{x}$ is a positive equilibrium point and if $f(x)$ has only one change in sign 
from positive to negative as $x$ increases from zero, then model \eqref{single:general} 
is globally asymptotically stable.
\end{corollary}


\subsection{Unexploited single-species model}

The standard single-species logistic model for population $x$ at time $t$ is given by
\begin{equation}
\label{stand:log}
\dot{x} = rx\left(1 - \frac{x}{k}\right), 
\end{equation}
where $r$ denotes the growth rate and $k$ the carrying capacity, i.e., the maximum population 
that the environment can sustain indefinitely.

It is clear that \eqref{stand:log} is a particular case of \eqref{single:general}, where 
$$f(x) = r\left(1 - \frac{x}{k}\right).$$
This model has two equilibrium points $\tilde{x}_1 = 0$ and $\tilde{x}_2 = k$, where the 
latter is positive. Moreover, when $k > x > 0$, $f(x) > 0$. On the other hand, when 
$x > k$, $f(x) < 0$. Thus, by Theorem \ref{cond:single}, model \eqref{stand:log} is globally 
asymptotically stable. 


\subsection{Exploited single-species models}

An important problem in the management of an exploited population is the stability of the 
population. We can distinguish two types of constant harvesting policies: constant effort
policy and constant quota policy. In the latter, a constant number of animals in a population 
is harvested per unit of time. In a constant effort policy, the effort applied per unit 
of time to catch the animals in a population is fixed. These two types of constant harvesting 
policies usually have completely different outcomes regarding the qualitative behavior of a
system. Constant quota harvesting policies are usually unstable, while constant effort 
harvesting policies are usually stable \cite{book}. Now,
we show these results.


\subsubsection{Constant effort policy}

Let $e$ be a measure of the effort applied in harvesting a population. Usually we have 
$e = qg$, where $q$ is the catchability coefficient and $g$ is the number of boats, men or 
other units used in the harvesting of the population. In these terms, the exploited population
is defined as
\begin{equation}
\label{const:effort}
\dot{x} = xf(x) - ex.
\end{equation}
As already proved and stated in Theorem~\ref{cond:single}, if $f(x)$ is strictly monotonic 
decreasing as $x$ increases and $\dot{x} = xf(x)$ has a positive equilibrium, then the 
unexploited population is globally asymptotically stable. In this case, we have that 
$\frac{\partial f}{\partial x} < 0$ for $x > 0$. 
Clearly, if  this is satisfied, then
$$\frac{\partial(f-e)}{\partial x} < 0,\, \text{for}\,\, x > 0.$$
Thus, if $f(x)$ is strictly monotonic decreasing for $x > 0$, the exploited and unexploited 
populations have the same qualitative behavior, provided in each case there is a positive 
equilibrium.

Let us prove that these conclusions apply to the logistic growth model considered so far. 
The logistic model with constant effort policy is given by
\begin{equation}
\label{log:effort}
\dot{x} = rx\left(1 - \frac{x}{k}\right) - ex = x\left[r\left(1 - \frac{x}{k}\right)
- e\right].
\end{equation}
Model \eqref{log:effort} has two equilibrium points $\tilde{x}_1 = 0$ and 
$\tilde{x}_2 = k\left(1 - \frac{e}{r}\right)$, which is positive if $e < r$. The latter comes 
from
$$r\left(1 - \frac{x}{k}\right) - e = 0.$$
Consider 
$$
g(x) = r\left(1 - \frac{x}{k}\right) - e,
$$ 
which is a linear function with a negative slope. So it is clear that when $\tilde{x}_2 > x > 0$, 
$g(x) > 0$, and when $x > \tilde{x}_2$, $g(x) < 0$. Thus, we are in conditions of 
Theorem~\ref{cond:single}, and so model \eqref{log:effort} is globally asymptotically stable.	 


\subsubsection{Constant quota policy}

Now, let us assume that a constant quota policy 
is prescribed for the harvesting of a 
single-species population. Let $h$ be the constant rate of harvest. In this case, a model 
of the harvested population is
\begin{equation}
\label{const:quota}
\dot{x} = xf(x) - h.
\end{equation}

For simplicity, let us consider the logistic model with constant quota policy. In this case, 
we have
\begin{equation}
\label{log:quota}
\dot{x} = rx\left(1 - \frac{x}{k}\right) - h. 
\end{equation}

Calculating the equilibrium points of \eqref{log:quota}, we get
$$ x = \frac{1}{2}\left(k\, \pm\, \sqrt{k^2 - \frac{4hk}{r}} \right).$$
If $h > \frac{rk}{4}$, the equilibrium points are not defined in $\mathbb{R}$, and thus we 
discard this case. To further study the remaining cases, first note that \eqref{log:quota}
can be rewritten as 
$$\dot{x} = x\left(r\left(1 - \frac{x}{k}\right) - \frac{h}{x}\right).$$

Let $f(x) = f_1(x) - f_2(x)$, where $f_1(x) = r\left(1 - \frac{x}{k}\right)$ and 
$f_2(x) = \frac{h}{x}$. An effective way to establish the stability of such exploited model is 
to plot both $f_1(x)$ and $f_2(x)$ on the same graphic. 

\begin{figure}[ht]
\centering
\includegraphics[scale=0.7]{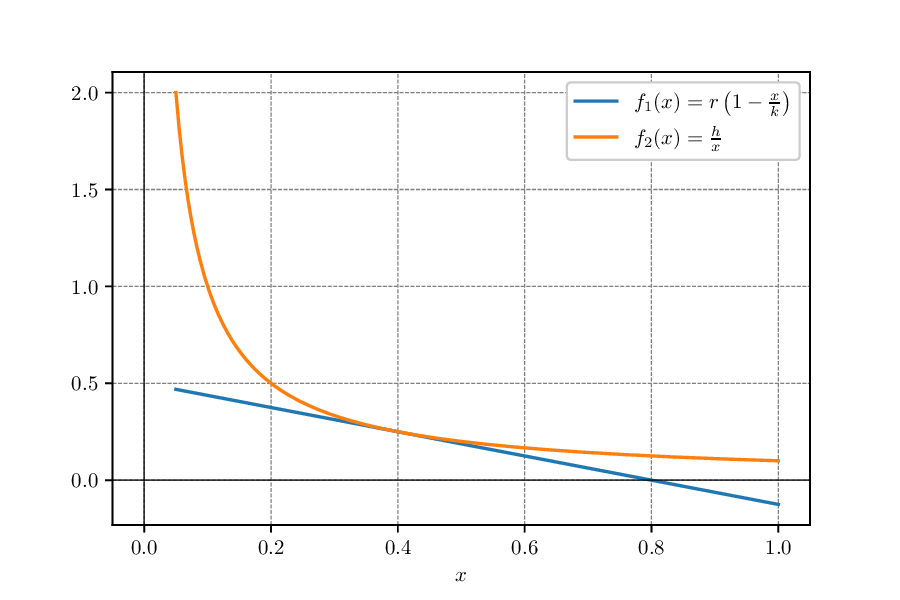}
\caption{A constant quota harvesting policy that generates one single equilibrium point.
Illustration with $h=0.1$, $r=0.5$, and $k=0.8$.}
\label{fig2}
\end{figure}

First, consider the case where $h = \frac{rk}{4}$. There is a single equilibrium point at 
$\tilde{x} = \frac{k}{2}$. For $h = 0.1$, $r = 0.5$ and $k = 0.8$, the behavior of $f_1(x)$ and 
$f_2(x)$ can be seen in Figure \ref{fig2}. It is clear that, for all $x$ except at the 
equilibrium $x = \frac{k}{2}$, $f_2(x)$ is always greater than $f_1(x)$. This means that 
$f(x)$ is negative for all $x$, with $x \neq \tilde{x}$, which means that for $0 < x < \tilde{x}$,
$\dot{V}(x) > 0$. Thus, following Theorem \ref{thm:unstable}, $\tilde{x}$ is unstable. 

\begin{figure}[ht]
\centering
\includegraphics[scale=0.7]{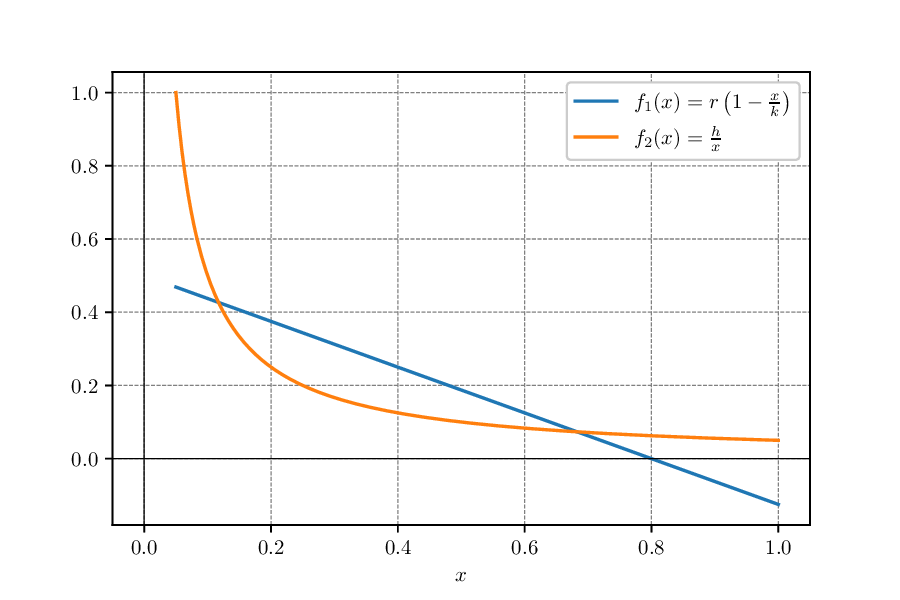}
\caption{A constant quota harvesting policy that generates two equilibria.
Illustration with $h=0.05$, $r=0.5$, and $k=0.8$.}
\label{fig3}
\end{figure}

On the other hand, when $h < \frac{rk}{4}$, there are two positive equilibria:
$$\tilde{x}_1 = \frac{1}{2}\left(k - \sqrt{k^2 - \frac{4hk}{r}} \right), \quad
\tilde{x}_2 = \frac{1}{2}\left(k + \sqrt{k^2 - \frac{4hk}{r}} \right).$$
As an example, for $h = 0.05$, $r = 0.5$ and $k = 0.8$, the behavior of $f_1(x)$ and 
$f_2(x)$ is sketched in Figure \ref{fig3}. We can see that
\begin{itemize}
\item for $0 < x < \tilde{x}_1$, $f_1(x) < f_2(x) \Rightarrow f(x) < 0$;
\item for $\tilde{x}_1 < x < \tilde{x}_2$, $f_1(x) > f_2(x) \Rightarrow f(x) > 0$;
\item for $x > \tilde{x}_2$, $f_1(x) < f_2(x) \Rightarrow f(x) < 0$. 
\end{itemize}
Thus, following the same reasoning as before, $\tilde{x}_1$ is unstable (Theorem 
\ref{thm:unstable}). Moreover, the region of stability for the equilibrium at 
$\tilde{x}_2$ is $\left\{x : x > \tilde{x}_1\right\}$.

\section{Optimal control of the logistic equation} 
\label{sec3}

The most popular way to model the behavior of an ecosystem is by means of a set of 
differential equations. It follows that the appropriate tool for formulating optimal 
policies for such systems is (the continuous time) optimal control theory. 

Generally, an optimal control problem consists on minimizing or maximizing 
a cost functional subject to a system and maybe to some initial or/and final 
conditions. Let $x$ be the state vector. In an ecosystem the components of 
$x(t) = (x_1(t), x_2(t), \dots, x_n(t))$ can be population densities at each instant of
time $t \in [a,b]$, with $a \geq 0$. To manage an ecosystem we must be able to manipulate
some variables which affect the dynamics of an ecosystem. Such variables are called 
controls. Let $u(t) = (u_1(t), u_2(t), \dots, u_m(t))$ be the control vector, $m \leq n$. 
In most applications, the components of $u(t)$ must satisfy certain constraints. 
Typically, one has $u(t) \in U$ with
\begin{equation}
\label{eq:def:U}
U = \{u \in \mathbb{R}^m:\, \alpha_i \leq u_i \leq \beta_i,\,i = 1, 2, \dots, m\}. 
\end{equation}

We now define an optimal control problem with given initial and terminal states on a
fixed finite time interval \cite{pontryagin}.

\begin{definition}
An optimal control problem with given initial and terminal states on a
fixed finite time interval $[a,b]$ consists in
\begin{equation}
\label{eq:cost}
\max\, \mathcal{J}[x(\cdot), u(\cdot)]= \int_{a}^{b} L(t,x(t),u(t))\, dt,
\end{equation}
subject to the differential system
\begin{equation}
\dot{x}(t) = f(t, x(t), u(t)),
\end{equation}
with boundary conditions
\begin{equation}
\label{boundary:cond}
x(a) = x_a\quad \text{and}\quad x(b) \geq x_b,
\end{equation}
and with control values
\begin{equation}
\label{eq:cont:values}
u(t) \in U \subset \mathbb{R}^m.
\end{equation}
\end{definition}
Condition \eqref{boundary:cond} refers to an initial condition, which simply means we can
measure the state vector at initial time, and a target $x(b) \geq x_b$. The specification 
of a target in an optimal control problem is itself a decision problem. But here, for 
simplicity, we assume that the target $x_b$ is specified in a given optimal control problem.  

In the following theorem, we recall a well-known necessary optimality condition: 
Pontryagin's Maximum Principle (PMP) \cite{pontryagin}. But before stating such result,
let us define the Hamiltonian associated with the optimal control problem 
\eqref{eq:cost}--\eqref{eq:cont:values} as
\begin{equation}
\label{hamiltonian}
H(t,x,u,\lambda_0,\lambda) = \lambda_0 L(t,x,u) + \lambda f(t,x,u).
\end{equation} 
We can either have $\lambda_0 = 0$ or $\lambda_0 = 1$. The functions $\lambda(t) = 
(\lambda_1(t), \lambda_2(t), \dots, \lambda_n(t))$ are called the adjoint functions. If 
$\lambda_0 = 0$, the corresponding optimal set $(x^*, u^*)$ is said to be abnormal.  
If $\lambda_0 = 1$, the optimal set $(x^*, u^*)$ is said to be normal.

\begin{theorem}[See \cite{pontryagin}]
\label{thm:pont}
A necessary condition for an admissible set $(x^*, u^*)$ to be optimal is that there exists 
$(\lambda_0, \lambda) \neq 0$ that satisfy
\begin{enumerate}
\item the adjoint system
$$
\dot{\lambda}(t) = -\frac{\partial H}{\partial x}(t,x^*(t),u^*(t), \lambda_0, \lambda(t));
$$
\item the optimality condition
$$
u^*(t) = \max_{u \in U} H(t, x^*(t), u, \lambda_0, \lambda(t));
$$  
\item the transversality condition
\begin{itemize}
\item $\lambda(b) = 0$ if $x(b) > x_b$;
\item $\lambda(b) = k$, with $k \in \mathbb{R}$, if $x(b) = x_b$;
\end{itemize}
\item and if a control variable switches at time $t$, then
$$\lambda(t-) = \lambda(t+),$$
$$
H(t-,x(t-), u(t-), \lambda_0, \lambda(t-)) = H(t+,x(t+), u(t+), \lambda_0, \lambda(t+)),
$$
where $\lambda(t-) = \lim_{\varepsilon \rightarrow 0^+} \lambda(t-\varepsilon)$
and $\lambda(t+) = \lim_{\varepsilon \rightarrow 0^+} \lambda(t+\varepsilon)$.
\end{enumerate} 
\end{theorem}	

\begin{definition}
The quadruple $(x^*(t), u^*(t), \lambda_0, \lambda(t))$ that satisfies 
Theorem~\ref{thm:pont} is called an extremal (or Pontryagin extremal).
\end{definition}

\begin{corollary}
\label{cor:opt}
If $U$ is defined as in \eqref{eq:def:U}, the optimality 
condition of Theorem~\ref{thm:pont} implies that
\begin{itemize}
\item $u^*_i(t) = \alpha_i$, only if 
$$
\frac{\partial H}{\partial u_i}(t, x^*(t), u^*(t), \lambda_0, \lambda(t)) < 0;
$$
\item $u^*_i(t) = \beta_i$, only if 
$$
\frac{\partial H}{\partial u_i}(t, x^*(t), u^*(t), \lambda_0, \lambda(t)) > 0;
$$
\item $\alpha_i < u^*_i(t) < \beta_i$, only if 
$$
\frac{\partial H}{\partial u_i}(t, x^*(t), u^*(t), \lambda_0, \lambda(t)) = 0
$$
and 
\begin{equation}
\label{matrix:sing}
\left(\frac{\partial^2 H}{\partial u_i\partial u_j}(t, x^*(t), u^*(t), \lambda_0, \lambda(t))\right)
\end{equation}
is a negative semidefinite matrix. 
\end{itemize}
\end{corollary}

If \eqref{matrix:sing} is a singular matrix, then $u^*$ is said to be a singular control. This 
usually occurs when one or more controls appears linearly both in the system dynamics and the 
objective functional. The most useful set of necessary conditions for singular control consists 
of the generalized Legendre condition, whose general form was obtained by Goh in 
\cite{goh:singular}. 

For simplicity, consider the one-dimensional case ($m = 1$)
and let $u(\cdot)$ be a single control that is associated 
with an Hamiltonian function $H(t, x, u, \lambda_0, \lambda)$ that is linear in $u$. 
The following theorem is called the generalized Legendre or Goh's conditions.

\begin{theorem}[See \cite{book}]
\label{thm:legendre}
If along an extremal one has
$$
\frac{\partial^2H}{\partial u^2} = 0,
$$
then necessary conditions for optimality are 
\begin{enumerate}
\item 
\begin{equation}
\label{eq:Goh:01}
\frac{\partial}{\partial u}
\left(\frac{d}{dt}
\left(\frac{\partial H}{\partial u}\right)\right) = 0\,\, \text{for all}\,\, 
t \in [a_1,b_1] \subset [a,b];
\end{equation}
\item if \eqref{eq:Goh:01} is satisfied, then 
$$-\frac{\partial}{\partial u}
\left(\frac{d^2}{dt^2}
\left(\frac{\partial H}{\partial u}\right)\right)$$
must be nonpositive for all $t \in [a_1,b_1] \subset [a,b]$.
\end{enumerate}
\end{theorem}

Here, we are interested to examine the use of optimal control theory to formulate optimal policies for the 
exploitation of a population described by the standard logistic model. More precisely, 
our optimal control problem is defined in the interval $[0,b]$ as follows:
\begin{equation}
\label{control:log}
\begin{aligned}
\max\,\, &\mathcal{J}[u(\cdot)] = \int_{0}^{b} u(t)\, dt,\\
&\dot{x}(t) = rx(t) \left(1 - \frac{x(t)}{k}\right) - u(t),\\
&x(0) = x_0,\\
&x(b) \geq x_b, \\
& 0 \leq u(t) \leq u_{\max},
\end{aligned}
\end{equation}
where $x_0, x_b \in \mathbb{R}$ and $u_{\max} \in \mathbb{R}$ are given.

By definition, the Hamiltonian is defined as
$$
H(x,u,\lambda_0, \lambda) = 
\lambda_0 u + \lambda\left(rx \left(1 - \frac{x}{k}\right) - u\right).
$$ 

Thus, the adjoint equation is given by
\begin{equation}
\label{adjoint:eq}
\dot{\lambda} = - \frac{\partial H}{\partial x}(x,u,\lambda_0,\lambda) 
= -\frac{\lambda r}{k}(k - 2x).
\end{equation}

Following Corollary~\ref{cor:opt}, the optimality conditions are as follows:
\begin{itemize}
\item $u = u_{\max}$ only if $\frac{\partial H}{\partial u} > 0$;
\item $u = 0$ only if $\frac{\partial H}{\partial u} < 0$;
\item $0 < u < u_{\max}$ only if $\frac{\partial H}{\partial u} = 0$.
\end{itemize}

Assuming $\lambda_0 = 1$, we have
$$H(x,u,\lambda_0,\lambda) = u + \lambda \left(rx\left(1 - \frac{x}{k}\right) - u\right),$$
which leads to
\begin{equation}
\label{opt:condition}
\frac{\partial H}{\partial u} = 1 - \lambda.
\end{equation} 

The terminal condition $x(b) \geq x_b$ is equivalent to the condition $x(b) = x_b$  or 
$x(b) > x_b$. Now, we look at these two different cases separately.

\begin{enumerate}
\item Let us start by analyzing the case $x(b) > x_b$. Here, the transversality condition 
gives $\lambda(b) = 0$. The adjoint equation is a first order homogeneous equation in 
$\lambda(t)$. Thus, the terminal condition $x(b) > x_b$ implies $\lambda(t) = 0$, for all 
$t \in [0,b]$. In fact, the solution of \eqref{adjoint:eq} is of the form
$$\lambda(t) = \lambda(0)\,e^{-\frac{r}{k}\int (k - 2x)\, dt}.$$
So, the only way for $\lambda(b)$ to be zero, as it is by the transversality condition, is 
if $\lambda(t) = 0$, for all $t \in [0,b]$. In this case, \eqref{opt:condition} simply 
becomes $$\frac{\partial H}{\partial u} = 1 > 0.$$
Thus, following the optimality conditions stated above, $u^* = u_{\max}$ in this case. 

\item Let us now suppose we have $x(b) = x_b$. In this case, the transversality condition 
becomes $\lambda(b) = c$, $c\in \mathbb{R}$, which is trivial and adds no information 
to our problem. Note that the control variable $u$ appears linearly in the Hamiltonian
function $H(x, u, \lambda_0, \lambda)$. This means that this particular problem may have a singular
control. 

First note that
$$\frac{\partial^2 H}{\partial u^2} = 0.$$ 
Along a singular extremal, we know that \eqref{opt:condition} is equal to zero. Thus, 
\begin{equation}
\label{first:deriv}
\frac{d}{dt}\left(\frac{\partial H}{\partial u}\right) = -\dot{\lambda}
= \frac{\lambda r}{k}(k - 2x) = 0,
\end{equation}
and 
\begin{equation}
\label{second:deriv}
\begin{aligned}
\frac{d^2}{dt^2}\left(\frac{\partial H}{\partial u}\right) 
&= \frac{\dot{\lambda}r}{k}(k-2x) + \frac{\lambda r}{k}(-2\dot{x})\\
&= -\frac{r}{k}(k-2x)^2 - \frac{2rx}{k}(k-x) + 2u = 0.
\end{aligned}
\end{equation}
Finally, 
\begin{equation}
\label{final:deriv}
-\frac{\partial}{\partial u}\left(\frac{d^2}{dt^2}
\left(\frac{\partial H}{\partial u}\right)\right) = -2 < 0.
\end{equation} 
So, we are in conditions of Theorem \ref{thm:legendre}, i.e., Goh's 
conditions for singular control are satisfied. 

We know that the optimal control maximizes $H(x,u,\lambda_0,\lambda)$
as a function of $u$. Thus, from \eqref{first:deriv}, we get the singular extremal
$$
x = \frac{k}{2}.
$$
Moreover, with $x = \frac{k}{2}$, from \eqref{second:deriv} we obtain the singular control
$$
u = \frac{rk}{4}.
$$
\end{enumerate}

Now, we need to make use of these conditions to define the optimal control policy along time.
As we shall see, a no switch policy is only feasible for some initial states. 

First, let us suppose, that $u_{\max} > \frac{rk}{4}$. Then, the obtained pair 
$$\left(x = \frac{k}{2},\, u = \frac{rk}{4}\right)$$ 
is admissible. 

For initial states above $\frac{k}{2}$, the singular control can be applied 
from the beginning since, after some time, it will lead to the singular extremal. 
On the other hand, for initial states below the singular extremal, it is needed to first 
consider $u = 0$, until the population density reaches $\frac{k}{2}$, at which 
point the singular control is applied, keeping the population constant from that time.
Note that these results are in line with the ones obtained in Section~\ref{sec2}, when
we analyzed the stability of the logistic model with constant quota policy. There, we 
concluded that when we had a single equilibrium $\tilde{x} = \frac{k}{2}$, for 
$0 < x < \tilde{x}$, $\dot{V}(x)$ was positive, leading to the instability of such 
equilibrium. In this case, any perturbation which causes the population to fall bellow 
$\frac{k}{2}$ will lead quickly to the extinction of the population unless a reduced control 
policy is imposed to allow the population to rise to the $\frac{k}{2}$ level.

\begin{figure}[ht]
\centering
\includegraphics[scale=0.7]{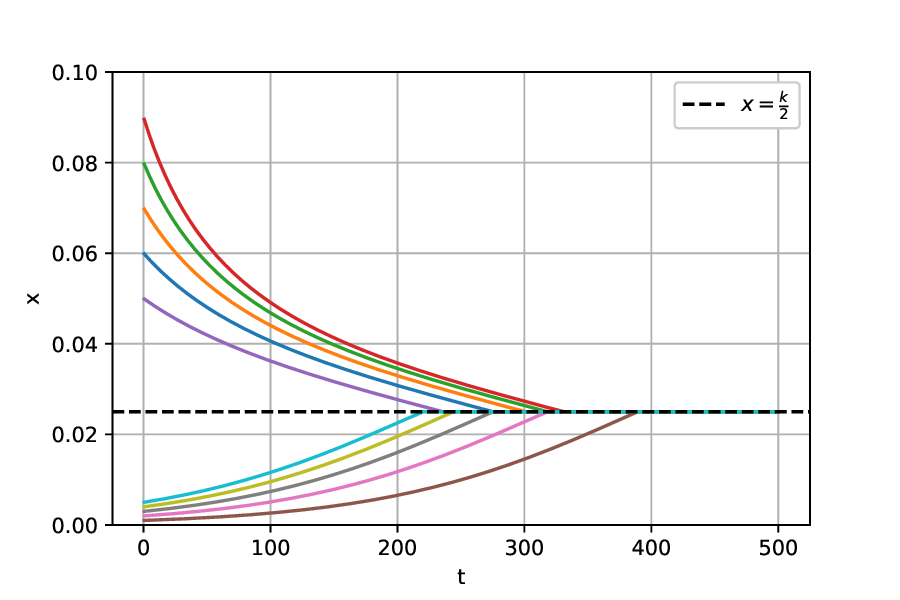}
\caption{Optimal trajectories for problem \eqref{control:log}	
when $u_{\max} > \frac{rk}{4}$.
Illustration with $b=500$, $x_b=0.001$, $r=0.01$, $k=0.05$, and $u_{\max} = 0.0002$.}
\label{fig4}
\end{figure}

Clearly, the functional is maximized if it is possible to use the maximum control
for all $t \in [0,b]$. Moreover, the terminal condition $x(b) > x_b$ implies that the optimal 
control $u = u_{\max}$ without any switch in the control variable. But, keeping the control 
equal to $u_{\max}$ for all time, is only feasible for some initial states. Thus, for 
initial states starting above the singular extremal, the simplest way to get to that extremal 
is to employ the control $u = u_{\max}$ until the trajectory intersects the curve 
$x = \frac{k}{2}$. When this intersection is reached, the control is changed to 
$u = \frac{rk}{4}$. 
The according behavior can be seen in Figure \ref{fig4}. 

\begin{figure}[ht]
\centering
\includegraphics[scale=0.7]{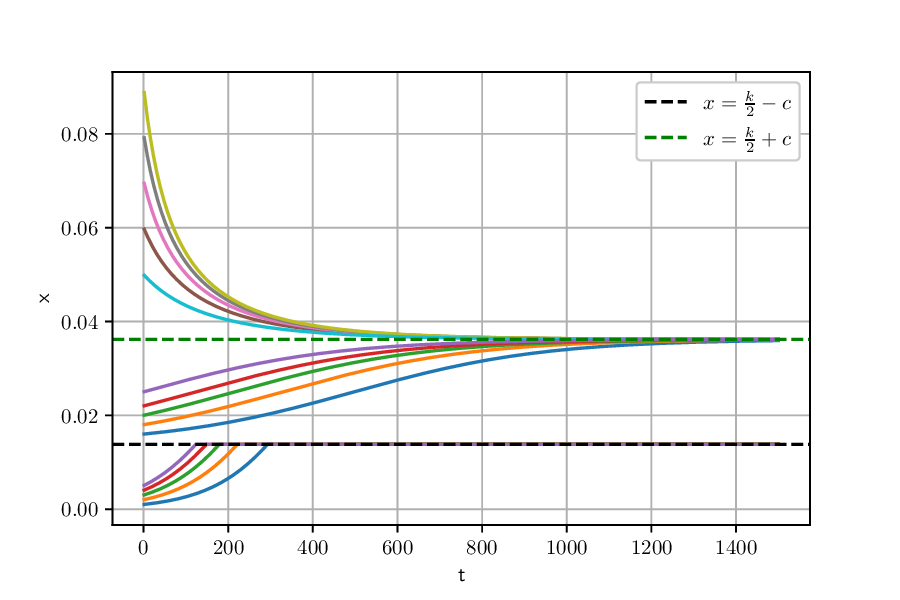}
\caption{Optimal trajectories for problem \eqref{control:log}
when $u_{\max} < \frac{rk}{4}$. Illustration with $b=1500$, $x_b=0.001$,
$r=0.01$, $k=0.05$, and $u_{\max} = 0.0001$.}
\label{fig5}
\end{figure}

Now, suppose that $u_{\max} < \frac{rk}{4}$. This way the singular control violates the
condition $0 < u < u_{\max}$, and so the pair 
$$\left(x = \frac{k}{2},\, u = \frac{rk}{4}\right)$$ 
is no longer admissible. Relating these results back to those obtained in 
Section~\ref{sec2}, when harvesting was lower than $\frac{rk}{4}$, we obtained two positive equilibrium 
points. The same happens here for $u < \frac{rk}{4}$ for all $t$. Such equilibria are equal to
$$ x = \frac{1}{2}\left(k\, \pm\, \sqrt{k^2 - \frac{4hk}{r}} \right).$$
A simple and stable policy in this case is to apply $u = u_{\max}$ whenever 
$$x > \frac{1}{2}\left(k\ - \sqrt{k^2 - \frac{4hk}{r}} \right)$$ 
and $u = 0$ whenever 
$$x < \frac{1}{2}\left(k\ - \sqrt{k^2 - \frac{4hk}{r}} \right),$$
point at which we change the control to $u = u_{\max}$. This policy is sketched in 
Figure~\ref{fig5}. 

\section{Logistic equation on time scales}
\label{sec4}

The theory of time scales was introduced in 1988
by Stefan Hilger \cite{HilgerPhD}. 
Such theory unifies continuous and discrete time 
into a single more general framework \cite{time:scales}. 
Recently, it has been applied to several models 
in biology and ecology
\cite{MR3636114,MR3895105,MR4738273}.

A time scale $\mathbb{T}$ is any nonempty closed subset of the real numbers. 
The forward jump operator is denoted by $\sigma(t)$ and $x^\Delta(t)$ is called 
the delta (or Hilger) derivative. The notation $x^\sigma = x \circ \sigma$ is often used.
If $\mathbb{T} = \mathbb{R}$, then $x^\sigma(t) = x(t)$ 
and the delta derivative is consistent with the classical derivative, 
i.e., $x^\Delta(t) = x'(t)$, for 
$t \in \mathbb{T} = \mathbb{R}$. If $\mathbb{T} = \mathbb{Z}$, 
then $x^\sigma(t) = x(t+1)$ and the delta derivative becomes 
the discrete analogue of a derivative, that is, $x^\Delta(t) = x(t+1) - x(t)$.
The reader interested in the calculus on time scales 
is referred to the book \cite{time:scales}.

As far as we are aware, there have been limited positive advancements in constructing the 
logistic equation on time scales \cite{sabrina}. 
Although there have been some studies about this dynamic equation, 
the authors have not succeed in ensuring its consistency. Precisely, let us consider 
Streipert's recent work \cite{sabrina}, where the following
analogue of the logistic growth model is given:
\begin{equation}
\label{time:scales:log}
x^\Delta = rx^\sigma\left(1 - \frac{x}{k}\right), 
\end{equation}
$x(t_0) = x_0$, with growth rate $r > 0$, 
carrying capacity $k > 0$, and initial population size 
$x_0 > 0$ at time $t_0 \in \mathbb{T}$. In the particular case
$\mathbb{T} = \mathbb{Z}$, one obtains
the discrete version of \eqref{time:scales:log} as
\begin{equation}
\label{discrete:log}
x(t+1) - x(t) = rx(t+1)\left(1 - \frac{x(t)}{k}\right),
\end{equation}
which is equivalent to 
$$
x(t+1) = \frac{x(t)}{1 - r\left(1 - \frac{x(t)}{k}\right)}.
$$
However, this discrete-time model \eqref{discrete:log} is not consistent with the classical continuous
logistic model: the solution can be negative, which does not make sense in the applications.
As an example, consider $x(0) = 2$, $r = 2$ and $k = 5$. This way, we get $x(1) = -10$.
This means that even when we start with an initial population size $x(0) > 0$, the model 
\eqref{discrete:log} fails to guarantee the non-negativity of solutions, which should be inherent to 
any system that intends to model population densities. This also shows that the time-scale
model \eqref{time:scales:log} of \cite{sabrina} is inconsistent.
To correct this inconsistency, we propose a different time-scale logistic model. 
For that we follow the rules stated by Mickens in \cite{mickens}. This method is being successfully applied
nowadays, since it generates discrete-time systems whose properties and qualitative behavior 
are identical to their continuous counterparts, as desired \cite{LemosSilva2023,MR4539989}.

The logistic time scales model we propose consists 
in changing the position of the composition of the state variable $x$
with the forward jump operator $\sigma$ as follows:
\begin{equation}
\label{time:scales:log:our}
x^\Delta = rx\left(1 - \frac{x^\sigma}{k}\right), \quad x(t_0) = x_0. 
\end{equation}
The study of the time-scale logistic model \eqref{time:scales:log:our} is not simple
and is under investigation. Here we just note that the corresponding discrete-time model,
obtained choosing  $\mathbb{T} = \mathbb{Z}$, is already consistent with the standard
logistic model:
\begin{equation}
\label{discrete:log:our}
x(t+1) - x(t) = rx(t)\left(1 - \frac{x(t+1)}{k}\right),
\end{equation}
which is equivalent to 
\begin{equation}
\label{discrete:log:our2}
x(t+1) = \frac{(r+1) k x(t)}{k+ r x(t)}.
\end{equation}
Clearly, beginning with $x(t_0) > 0$, now with \eqref{discrete:log:our2}
one has $x(t) > 0$ for all $t > 0$.
 
 
\section*{Acknowledgements}

The authors were partially supported by 
the Portuguese Foundation for Science and Technology (FCT):
Lemos-Silva and Torres through
the Center for Research and Development in Mathematics 
and Applications (CIDMA), project UIDB/04106/2020 
(\url{https://doi.org/10.54499/UIDB/04106/2020});
Vaz through the Center of Mathematics and Applications 
of \emph{Universidade da Beira Interior} (CMA-UBI), 
project UIDB/00212/2020 (\url{https://doi.org/10.54499/UIDB/00212/2020}). 
Lemos-Silva is also supported by the FCT PhD fellowship
UI/BD/154853/2023; Torres within the project 
``Mathematical Modelling of Multiscale Control Systems: 
Applications to Human Diseases'' (CoSysM3), 
Reference 2022.03091.PTDC (\url{https://doi.org/10.54499/2022.03091.PTDC}), 
financially supported by national funds (OE) through FCT/MCTES.
The authors are very grateful to a referee for the careful reading of the paper and for
several comments and detailed suggestions that helped to improve 
the exposition.



\end{document}